\documentclass{lncse}
\font\elevenbb=msbm10 at 10.95pt

\def\K{\hbox{\elevenbb K}}
\def\N{\hbox{\elevenbb N}}
\def\M{\hbox{\elevenbb M}}
\def\R{\hbox{\elevenbb R}}

\def\Gr{Gr\"obner }
\def \bg #1 {\begin{tabular}{{#1}}}
\def \nd {\end{tabular}}
\newcommand \mwhile {{\bf while}\hspace{0.25cm}}
\newcommand \mrepeat {{\bf repeat}\hspace{0.25cm}}
\newcommand \muntil {{\bf until}\hspace{0.25cm}}
\newcommand \mfore {{\bf for\hspace{0.25cm}each}\hspace{0.25cm}}
\newcommand \mdo {{\bf do}\hspace{0.25cm}}

\newcommand \mif {{\bf if}\hspace{0.25cm}}
\newcommand \mthen {{\bf then}\hspace{0.25cm}}
\newcommand \melse {{\bf else}\hspace{0.25cm}}
\newcommand \mchoose {{\bf choose}\hspace{0.25cm}}
\newcommand \mand {{\bf and}\hspace{0.25cm}}

\newcommand \mbegin {{\bf begin}}
\newcommand \mend {{\bf end}}
\newcommand \bb {\hspace{0.25cm}}
\newcommand \h {\hspace{0.25cm}}
\newcommand \hh {\hspace{0.5cm}}
\newcommand \hhh {\hspace{0.75cm}}
\newcommand \hhhh {\hspace{1.0cm}}
\newcommand \hhhhh {\hspace{1.25cm}}
\newcommand \hhhhhh {\hspace{1.5cm}}
\newcommand \hhhhhhh {\hspace{1.75cm}}
\newcommand \hhhhhhhh {\hspace{2.0cm}}

\newcounter{cc}
\newcommand \hln {\hfill \addtocounter{cc}{1} \arabic{cc}
            \vskip 0.0cm \noindent }
\begin{document}
\title{Completion of Linear Differential Systems to
  Involution\thanks{
  {\em
  Computer Algebra in Scientific
  Computing}, V.G.Ganzha, E.W.Mayr and E.V. Vorozhtsov (Eds.),
  Springer-Verlag, Berlin, 1999, pp.115-137.
}}
\author{Vladimir P.
Gerdt} \institute{Laboratory of Computing Techniques and Automation,
Joint Institute for Nuclear Research, 141980 Dubna, Russia}
\maketitle
\begin{abstract}
In this paper we generalize the involutive methods and algorithms
devised for polynomial ideals to differential ones
generated by a finite set of linear differential polynomials in the
differential polynomial ring over a zero characteristic differential
field. Given a ranking of derivative terms and an involutive division,
we formulate the involutivity conditions which form a basis of
involutive algorithms. We present an algorithm for computation of
a minimal involutive differential basis. Its correctness and termination
hold for any constructive and noetherian involutive division. As two
important applications we consider posing of an initial value problem
for a linear differen\-tial system providing uniqueness of its
solution and the Lie symmetry analysis of nonlinear differential
equations. In particular, this allows to determine the structure of
arbitrariness in general solution of linear systems and thereby to
find the size of symmetry group.
\end{abstract}
\section{Introduction}
Among the properties of systems of analytical partial differential equations
(PDEs) which may be
investigated without their explicit integration there are
compatibility and formulation of an {\em initial-value problem}
providing existence and uniqueness of the solution. The
classical Cauchy-Kowalevsky theorem establishes a certain class of
quasilinear PDEs which admit posing such an initial-value problem. The main
obstacle in investigating other classes of PDE systems of some given
order $q$ is existence of {\em integrability conditions}, that is,
such relations for derivatives of order $\leq q$ which are
differential but not pure algebraic consequences of equations in
the system.

An {\em involutive} system of PDEs has all the integrability
conditions incor\-porated in it. This means that prolongations of
the system do not reveal integrability conditions. Extension of a
system by its integrability conditions is called {\em completion}.
The concept of involutivity was invented hundred years ago
by E.Cartan~\cite{Cartan1} in his investigation of the Pfaff type
equations in total differentials. For these purposes he used the
exterior calculus developed by himself. The Cartan approach was
generalized by K\"{a}hler\cite{Kahler} to arbitrary systems of exterior
differential equations. The underlying completion
proce\-dure~\cite{Cartan} was implemented in~\cite{ASY,HT}.

In his study of the formal power series solutions of PDEs, Riquier
intro\-duced~\cite{Riquier} a class of relevant {\em rankings} for partial
derivatives and considered systems of {\em orthonomic} equations which are
solved with respect to the highest rank derivatives called
{\em principal}. Thereby, these derivatives, by the equations in the system,
are defined in terms of the other derivatives called {\em parametric}. An
integrability condition gives a constraint for parametric derivatives, and
that of them of the highest ranking becomes the principal derivative.
Recently Riquier's class of rankings was generalized in~\cite{RR97}.

Janet made the further development of Riquier's approach.
He observed \cite{Janet} that the integrability
conditions may occur only from prolongations with respect to
certain independent variables called {\em nonmultiplicative}.
Prolongations with respect to the rest of variables called
{\em multiplicative} never lead to integrability conditions. Given a set
of principal derivatives, Janet gave the prescription how to separate
variables into multiplicative and nonmultiplicative for every
equation in the system. He formulated, on this ground, the
{\em involutivity conditions} for orthonomic systems and designed an
algorithm for their completion. This approach to completion is known
as {\em Riquier-Janet theory} and was implemented
in~\cite{Schwarz1,Topunov,Schwarz2}.

A system satisfying the Janet involutivity conditions is often called
{\em passive}. This involutivity is generally coordinate dependent.
On the other hand, the modern {\em formal theory} of PDEs developed in
60s-70s by Spencer and others (see~\cite{Pommaret,Pommaret94}) allows to
formulate the involutivity intrinsically, in a coordinate
independent way. The formal theory relies on another definition of
multiplicative and nonmultiplicative variables which was known to
Janet as long ago as in 20s, but called nowadays after Pommaret
because of its importance in the technique presented
in~\cite{Pommaret}. The implementation in Axiom of completion
based on the formal theory was presented in~\cite{Seiler1,Seiler2}.

Thomas in~\cite{Thomas} used another separation of independent
variables into multiplicative and nonmultiplicative and generalized
the Riquier-Janet theory to non-orthonomic algebraic PDEs. Given a
system of PDEs, he showed that in a finite number of steps one can:
(i) check its compatibility; (ii) if the system is compatible, then
split it into a finite number of {\em simple} systems involving
generally both equations and inequalities and such that their
equation parts are orthonomic and can be completed to involution.
This splitting is similar to that generated by the
Rosenfeld-Gr\"{o}bner algorithm~\cite{Boulier}.

In paper~\cite{ZB93} for Pommaret separation of independent
variables it was shown that involutive (passive) basis of a
non-differential polynomial ideal is a \Gr basis. The
implementation in Reduce of the proposed completion algorithm for
polynomial bases demonstrated a high computational efficiency of the
involutive technique. However, Pommaret bases may not exist for
positive dimensional ideals unlike Janet and Thomas bases.

The above classical separations of variables into multiplicative
and nonmultiplicative are particular cases of {\em involutive monomial
division}, a concept invented and analyzed in~\cite{GB1} (cf.~\cite{Apel}).
The polynomial completion
algorithms designed for a general involutive division~\cite{GB1,GB2} were
implemented in Reduce for Pommaret division. Different
involutive divisions and completion of monomial
sets have also been implemented in Mathematica \cite{GBC}.
In \cite{G98} we generalized the algorithm of paper \cite{GB2} to
arbitrary completion ordering.

One more efficient method for the completion of linear PDEs to an involutive form
called {\em standard} which is not based on the separation of variables
was developed in~\cite{Reid1} and implemented in Maple. The
extension of this method to nonlinear PDEs is given in
paper~\cite{Reid3}.

In the present paper we generalize the involutive methods and algorithms
devised in~\cite{GB1,GB2,G98} for polynomial ideals to differential ideals
generated by a finite set of linear polynomials. We formulate the
involutivity conditions for the differential case. If a set satisfies
the involutivity conditions it is called an {\em involutive basis}.
Similar to the pure algebraic case, a linear involutive basis is a
differential \Gr basis~\cite{Carra,Ollivier} which is not generally reduced.
We present an algorithm for computation of a minimal involutive basis.
This algorithm is the straightforward generalization of the
polynomial involutive algorithm~\cite{GB2,G98}. As well as for the latter,
the correctness and termination of the former hold for any constructive and
noetherian involutive division.

An important application of the involutive method is posing
an initial value problem providing the unique solution of a system
of PDEs. For linear involutive systems we
formulate such an initial value problem and thereby generalize
the classical results of Janet~\cite{Janet} to arbitrary involutive
divisions. This formulation makes it possible, among other things, to
reveal the structure of arbitrariness in general solution. Given
a linear involutive basis, we write also the explicit formulae for the
Hilbert function and the Hilbert polynomial of the corresponding
differential ideal which are the straightforward generalizations of
their polynomial analogues \cite{Apel,GBC}.

Another important application of the new algorithm is the {\em Lie symmetry
analysis} of nonlinear differential equations. It is because of the fact
that completion to involution is the most general and universal method of
integrating the determining system of linear PDEs for infinitesimal Lie
symmetry generators~\cite{Hereman}. Moreover, an involutive form of
determining equations allows to construct the Lie symmetry algebra
without their explicit integration~\cite{Reid2}. In particular, for an
involutive determining system the size of symmetry group can easily be found
that was shown for Janet bases in~\cite{Schwarz2}. Though reduced \Gr
bases for the determining equations do not generally reveal information
on Lie symmetry groups, and more generally on the solution space,
so explicitly as involutive bases, they are also very useful for Lie
symmetry analysis as shown in~\cite{Mansfield}. The facilities of the
Maple package devised by the first author and used in the paper go far
beyond linear differential systems, and it can also be fruitfully
applied to nonlinear systems.

\section{Preliminaries}
Let $\R=\K\{y_1,\ldots,y_m\}$ be a differential polynomial
ring~\cite{Ritt,Kolchin} with the set of differential indeterminates
$\{y_1,\ldots,y_m\}$, and $\K\subset \R$
is a differential field of zero characteristic with a finite number
of mutually commuting derivation operators
$\partial /\partial x_1,\ldots,\partial /\partial x_n$.
Elements in $\R$ are differential polynomials in
$\{y_1,\ldots,y_m\}$. In this paper we use the following
notations and conventions:

$f,g,h,p\in \R$ are linear differential polynomials.

$F,G,H\subset \R$ are finite sets of linear differential polynomials.

${\cal F}=\{f=0\ |\ f\in F\}$  is a linear system of PDEs.

$\N=\{0,1,2,\ldots\}$ is the set of nonnegative integers.

$\alpha,\beta,\gamma \in \N^n$  are multiindices.

$lcm(\alpha,\beta)$ is the least common multiple of $alpha,beta$.

$X=\{x_1,\ldots,x_n\}$ is the set of independent variables.

$R=K[X]$ is the polynomial ring over the field $K$ of zero
characteristic.

$R\supset \M=\{x^\alpha=x_1^{\alpha_1}\cdots x_n^{\alpha_n}\ |\ \alpha_i\in \N\}$ is the set
of monomials in $X$.

$i=1,\ldots,n$ indexes derivation operators $\partial_i=\partial/\partial_{x_i}$.

$j=1,\ldots,m$ indexes indeterminates $y_j$.

$u,v,w$ are elements in $\M$.

$U,V\subset \M$ are finite monomial sets.

$(U)$ is the monomial ideal in $R$ generated by $U$.

$deg_i(u)$ is the degree of $x_i$ in $u\in \M$.

$deg(u)=\sum_{i=1}^n deg_i(u)$ is the total degree of $u$.

$\partial_{\alpha} y_j=\frac {\partial^{\alpha_1+\cdots+\alpha_n}}
{\partial x_1^{\alpha_1} \cdots \partial x_n^{\alpha_n}} y_j$
is a derivative.

$ord(\partial_\alpha y_j)=\sum_{i=1}^n \alpha_i$ is the order of $\partial_\alpha y_j$.

$\theta,\vartheta$ are derivatives.

$lcm(\partial_\alpha y_j,\partial_\beta y_j)=\partial_{lcm(\alpha,\beta)}y_j$.

$\prec $, $\prec_c$  are rankings of derivatives.

$ld(f)$ is the leading derivative in $f\in \R$.

$lc(f)\in K$ is the coefficient of $ld(f)$.

$ld(F)$ is the set of leading derivatives in $F\subset \R$.

$[F]$ is the differential ideal in $\R$ generated by $F$.

$L$ is an involutive division.

$L(u,U)$ is the set of ($L-$)multiplicative monomials for
 $u\in U$.

$NM_L(u,U)$ is the set of $L-$nonmultiplicative variables for $u\in U$.

$M_L(u,U)$ is the set of $L-$multiplicative variables for $u\in U$.

$x^\alpha\in \M$ is the monomial associated with the derivative
$\partial_\alpha y_j$.

$\cup_{j=1}^m U_j$ is the monomial set
associated with the set $ld(F)=\cup_{j=1}^m \{ld_j(F)\}$.

$W=\cup_{j=1}^m \{W_j\ |\ W_j\subset M\}$ is the complementary set for $\cup_{j=1}^m U_j$.

$\mathcal{G}$ is the set of $L-$generators of $W$.

$\vartheta = \partial_L \theta$ is a multiplicative prolongation
of $\theta$.

$\partial_{x_i} \cdot \theta$ is the nonmultiplicative prolongation of $\theta$
w.r.t. $x_i$.

$\partial_\alpha \cdot \theta$ is a nonmultiplicative prolongation of $\theta$.

$NF_L(p,F)$ is the $L-$normal form of $p$ modulo $F$.

$NM_L(f,F)\subseteq X$ is the set of
($L-$)nonmultiplicative variables for $f\in F$.

$C_L(F)=\cup_{\theta\in ld(F)} \{\vartheta\ |\
 \vartheta = \partial_L \theta\}$ is the $L-$cone generated by $F$.

\noindent
In this paper we distinguish two rankings (c.f.~\cite{G98}):
a {\em main ranking} and a {\em completion ranking} denoted by $\succ$ and $\succ_c$,
respectively. The main ranking will be used, as usually, for isolation of the
leading derivatives in differential polynomials whereas the completion ranking
serves for taking the lowest nonmultiplicative prolongations by the normal
strategy~\cite{GB1} and thereby controlling the property of partial
involutivity introduced in Sect.~4.

\section{Basic Concepts and Definitions}

Throughout this paper we exploit the well-known algorithmic similarities
between pure algebraic polynomial systems and linear differential
systems \cite{Pommaret94,Gerdt97}. In so doing, the basic algorithmic ideas
go back to Janet~\cite{Janet} who invented the constructive
approach to study of PDEs in terms of the corresponding monomial sets which
is based on the following association between derivatives and monomials:
\begin{equation}
\partial_\alpha y_j=\frac {\partial^{\alpha_1+\cdots+\alpha_n} y_j}
{\partial x_1^{\alpha_1} \cdots \partial x_n^{\alpha_n}} \Longleftrightarrow
x^\alpha=x_1^{\alpha_1}\cdots x_n^{\alpha_n}\,. \label{d-m}
\end{equation}
The monomials associated with the different indeterminates $y_j$ are to be considered
as belonging to different monomial sets $U_j\in \M$ indexed by
subscript $j$ of the indeterminate.

\begin{definition}~\cite{Kolchin} A total ordering $\prec$ over the set
of derivatives $\partial_\alpha y_j$ is called a {\em ranking} if
it satisfies: (i) $\partial_{i} {\partial_\alpha y_j} \succ {\partial_\alpha y_j}$,
(ii) ${\partial_\alpha y_j} \succ {\partial_\beta y_k} \Longleftrightarrow
{\partial_\gamma}{\partial_\alpha y_j}
   \succ {\partial_\gamma}{\partial_\beta y_k}$
for all $i,j,k,\alpha,\beta,\gamma$. A ranking $\prec$ is said to be
{\em orderly} if $\theta \succ \vartheta$ whenever
$ord(\theta)>ord(\vartheta)$.
\label{ranking}
\end{definition}

\noindent
The association~(\ref{d-m}) implies the reduction of a ranking $\prec$
to the associated admissible monomial ordering, and throughout the paper we shall
assume that

\begin{equation}
\partial_{1}\succ \partial_{2}\succ \cdots \succ \partial_{n} \Longleftrightarrow
x_1\succ x_2\succ\cdots\succ x_n\,. \label{var_order}
\end{equation}

\begin{remark} Given a finite set $F\subset \R$ and a ranking $\succ$,
 set $ld(F)$ of the leading derivatives is partitioned
 $ld(F)=\cup_jld_j(F)$ into subsets $ld_j(F)$ corresponding to different
 indeterminates $y_j$ which occur in $ld(F)$. For an involutive
 division $L$ defined as follows each subset generates for every its element the separation
 of independent variables into multiplicative and nonmultiplicative ones.
\label{rem_inv_subs}
\end{remark}

\begin{definition}~\cite{GB1}
An {\em involutive division} $L$ on $\M$ is given, if for any finite
monomial set $U\subset \M$ and for any $u\in U$ there is given a
submonoid $L(u,U)$ of $\M$ satisfying the conditions:
\begin{tabbing}
~~(a)~~\=If $w\in L(u,U)$ and $v|w$, then $v\in L(u,U)$. \\
~~(b)  \>If $u,v\in U$ and $uL(u,U)\cap vL(v,U)\not=\emptyset$, then
$u\in vL(v,U)$ \\
       \> or $v\in uL(u,U)$. \\
~~(c)  \> If $v\in U$ and $v\in uL(u,U)$, then $L(v,U)\subseteq L(u,U)$. \\
~~(d)  \> If $V\subseteq U$, then $L(u,U)\subseteq L(u,V)$ for all $u\in V$.
\end{tabbing}
Elements of $L(u,U)$ are called {\em multiplicative} for $u$.  If
$w\in uL(u,U)$, $u$  is called an involutive divisor or {\em $(L-)$divisor}
of $w$. In such an event the monomial $v=w/u$ is called
{\em $L-$multiplicative} for
$u$. If $u$ is a conventional divisor of $w$ but not $L-$divisor, then
$v$ is called {\em nonmultiplicative} for $u$.
\label{inv_div}
\end{definition}

\begin{remark} Definition~\ref{inv_div} for every $u\in U$ provides
the partition
\begin{equation}
X=M_L(u,U)
\cup NM_L(u,U),\quad M_L\subset L(u,U)
\label{partition}
\end{equation}
of the set of variables $X=\{x_1,\ldots,x_n\}$ into subset
$M_L(u,U)$ of {\em multiplicative variables} for $u$
and subset $NM_L(u,U)$ of the remaining
{\em non\-mul\-ti\-pli\-ca\-ti\-ve variables}.
Conversely, if for any finite set $U\subset \M$ and any $u\in U$ the
partition~(\ref{partition})  of variables into multiplicative and
nonmultiplicative is given such that the corresponding submonoid $L(u)$
satisfies the conditions (b)-(d) in Definition~\ref{inv_div}, then the
partition generates an involutive division.
\label{rem_sep}
\end{remark}

\begin{definition}~\cite{GB1} A monomial set $U$
 is called $L-${\em autoreduced} if
$uL(u,U)\cap vL(v,U)=\emptyset$ holds for all distinct $u,v\in U$.
\label{autored_ms}
\end{definition}

\begin{definition}~\cite{GB1} A monomial set $\tilde{U}$
 is called an $L-${\em completion} of a set $U\subseteq \tilde{U}$ if
$$
 (\forall u\in U)\ (\forall w\in \M)\
 (\exists v\in \tilde{U})\ \ [\ uw\in vL(v,\tilde{U})\ ]\,.
$$
 If there exists a finite $L-$completion $\tilde{U}$ of a finite set
 $U$, then the latter is called {\em finitely generated} with respect
 to $L$. The involutive division $L$ is {\em noetherian} if
 every finite set $U$ is finitely generated with respect to $L$.
 If $\tilde{U}=U$, then $U$ is called {\em $L-$complete}. An $L-$autoreduced
 and complete set is called {\em
 $(L-)$involutive}.
\label{id_noetherian}
\end{definition}

\begin{definition}~\cite{GB1} Given a monomial set $U$, the set
$\cup_{u\in U}\,u\,\M$ is called {\em the cone} generated by $U$ and
denoted by $C(U)$. The set $\cup_{u\in U}\,u\,L(u,U)$ is called
{\em the involutive cone} of $U$ with respect to $L$ and denoted by
$C_L(U)$.
\label{cone}
\end{definition}

\noindent
Thus, the set $\tilde{U}$ is an $L-$completion of $U$ if $C(\tilde{U})=C_L(\tilde{U})=C(U)$.
Correspondingly, for an involutive set $U$ the equality $C(U)=C_L(U)$ holds.

Whereas noetherity provides existence of a finite involutive basis for any
polynomial ideal, another important properties of an involutive division called
{\em continuity} and {\em constructivity} provide the algorithmic construction of
involutive bases~\cite{GB1}. Continuity implies involutivity when the local
involutivity holds whereas constructivity strengthens continuity and
allows to compute involutive bases by sequential examination of single
nonmultiplicative prolongations only. We refer to papers~\cite{GB1,GB2,G98}
for description of these topics in detail. In those papers some examples
of involutive divisions were studied (see also~\cite{GBCK}) which include
three divisions called after Janet, Thomas and Pommaret, because they have
used the corresponding separations of variables for involutivity analysis
of PDEs~\cite{Janet,Thomas,Pommaret}. Other two divisions called Division I
and II were introduced in~\cite{GB2}, and a class of involutive divisions
called Induced division, since every division in the class is induced
by an admissible monomial orderings, was introduced in~\cite{G98}. All
those divisions are constructive and, except Pommaret division, they are
noetherian. Below we use three of those divisions defined as follows.

\begin{definition} Janet division~\cite{Janet}.
Let $U\subset \M$ be a finite set. Divide $U$ into groups
labeled by non-negative integers $\alpha_1,\ldots,\alpha_i$
$(1\leq i\leq n)$ :
$$[\alpha_1,\ldots,\alpha_i]=\{\ u\ \in U\ |\ \alpha_j=deg_j(u),\
1\leq j\leq i\ \}.$$
Then $x_i$ is multiplicative for $u\in U$ if $i=1$
and $deg_1(u)=\max\{deg_1(v)\ |\ v\in U\}$, or $u\in
[\alpha_1,\ldots,\alpha_{i-1}]$ and $deg_i(u)=\max\{deg_i(v)\ |\ v\in
[\alpha_1,\ldots,\alpha_{i-1}]\}$ for $i>1$.
\label{div_J}
\end{definition}

\begin{definition}
 Pommaret division~\cite{Pommaret}. For a monomial
$u=x_1^{\alpha_1}\cdots x_k^{\alpha_k}$ with $\alpha_k>0$ the variables $x_j,j\geq k$ are
considered as multiplicative and the other variables as nonmultiplicative.
For $u=1$ all the variables are multiplicative.
\label{div_P}
\end{definition}

\begin{definition} Lexicographically induced division~\cite{G98}.
 A variable $x_i$ is nonmultiplicative for $u\in U$ if there is $v\in U$ such that
 $v\prec_{Lex} u$ and $deg_i(u)<deg_i(v)$, where $\succ_{Lex}$ denotes
 the lexicographical ordering.
\label{ind_div}
\end{definition}

\noindent
In the sequel Janet, Pommaret and Lexicographically induced divisions will be distinguished
by the subscripts
$J,P$ and $D_{Lex}$, respectively.

\begin{example} Separation of variables for
set $U=\{x_1^2x_3,x_1x_2,x_1x_3^2\}$ and
ordering~(\ref{var_order}) for the above defined
three divisions:

\begin{center}
\bg {|c|c|c|c|c|c|c|} \hline\hline
Element & \multicolumn{6}{c|}{Separation of variables}
\\ \cline{2-7}
in $U$ & \multicolumn{2}{c|}{Janet} &\multicolumn{2}{c|}{Pommaret} &
\multicolumn{2}{c|}{Lex. induced} \\ \cline{2-7}
 & $M_J$ & $NM_J$  & $M_P$   & $NM_P$   & $M_{D_{Lex}}$ & $NM_{D_{Lex}}$ \\ \hline
$x_1^2x_3$ & $x_1,x_2,x_3$ & $-$ & $x_3$ & $x_1,x_2$ & $x_1$ & $x_2,x_3$ \\
$x_1x_2$ & $x_2,x_3$ & $x_1$ & $x_2,x_3$ & $x_1$ & $x_1,x_2$ & $x_3$ \\
$x_1x_3^2$ & $x_3$ & $x_1,x_2$ & $x_3$ & $x_1,x_2$ & $x_1,x_2,x_3$ & $-$ \\
\hline \hline
\nd
\end{center}

\noindent
The corresponding $L$-completions of $U$ are
\begin{eqnarray*}
&& \tilde{U}_J=\{x_1^2x_3,x_1x_2,x_1x_3^2,x_1^2x_2\}, \\
&& \tilde{U}_P=\{x_1^2x_3,x_1x_2,x_1x_3^2,x_1^2x_2,\ldots,x_1^{i+2}x_2,
\ldots,x_1^{j+2}x_3,\ldots\}, \\
&& \tilde{U}_{D_{Lex}}=\{x_1^2x_3,x_1x_2,x_1x_3^2,x_1x_2x_3\}.
\end{eqnarray*}
where $i,j\in \N$. This example explicitly shows the non-noetherity of Pommaret
division.
\label{exm_1}
\end{example}

\begin{definition} Given a finite set $F\subset \R$, a ranking $\succ$ and
an involutive division $L$, the derivative $\vartheta = \partial_\beta y_j$
 will be called a {\em multiplicative prolongation} of $\theta = \partial_\alpha y_j\in ld_j(F)$
 and denoted by $\vartheta = \partial_L \theta$, if the associated monomials satisfy
 $x^\beta\in x^\alpha L(x^\alpha,U_j)$. Otherwise the prolongation will be called
 {\em nonmultiplicative}. Respectively, the corresponding prolongation
$\partial_{\beta}f$ of the element
 $f\in F$ with $ld(f)=\partial_\alpha y_j$ will be called
{\em multiplicative} and denoted by $\partial_L(f)$ or
{\em nonmultiplicative}.
 The set $C_L(F)=\cup_{\theta\in ld(F)} \{\vartheta\ |\ \vartheta = \partial_L \theta\}$
 will be called the {\em $L-$cone} generated by $F$. If $\partial_{i}f$ is
 a nonmultiplicative prolongation of $f\in F$, we shall write $x_i\in NM_L(f,F)$.
\label{def_prol}
\end{definition}

\section{Linear Involutive Differential Bases}

In this section we generalize the results obtained in
papers~\cite{GB1,G98} for commutative algebra to differential algebra of linear polynomials.
Proofs of
the theorems are omitted because of similarity with the proofs of
their algebraic analogues.

\begin{definition} Given an involutive division $L$, a
 finite set $F\subset \R$ of linear
 differential polynomials, a ranking $\succ$ and a linear polynomial $p\in \R$, we shall say:
 \begin{enumerate}
 \item $p$ is {\em $L-$reducible modulo} $f\in F$
  if $p$ has a term $a\,\theta$, $(a\in K\setminus \{0\})$
  such that $\theta=\partial_L ld(f)$. It yields the {\em $L-$reduction}
  $p\rightarrow g=p-\left(a/lc(f)\right)\partial_\beta f$ where $\partial_\beta ld(f)=\theta$.
 \item $p$ is {\em $L-$reducible modulo} $F$ if there is $f\in F$ such
  that $p$ is $L-$reducible modulo $f$.
 \item $p$ is {\em in $L-$normal form modulo $F$}, if $p$ is not $L-$reducible
 modulo $F$.
\end{enumerate}
\label{inv_red}
\end{definition}

\noindent
We denote the $L-$normal form of $p$ modulo $F$ by $NF_L(p,F)$.

As a $L-$normal form algorithm one can use the following differential analogue
of the polynomial normal form algorithm~\cite{GB1}:

\vskip 0.3cm
\noindent
\hh Algorithm {\bf InvolutiveNormalForm:}
\vskip 0.2cm
\noindent
\hh {\bf Input:}  $p,\,F,\,L,\,\prec$
\vskip 0.0cm \noindent
\hh {\bf Output:} $h=NF_L(p,F)$
\vskip 0.0cm \noindent
\hh \mbegin
\vskip 0.0cm \noindent
\hhh $h:=p$
\vskip 0.0cm \noindent
\hhh \mwhile exist $f\in F$ and a term $a\,\theta$ $\left(a\in \K\setminus \{0\}\right)$ of $h$
\vskip 0.0cm \noindent
\hhhh such that $\theta=\partial_L ld(f)$\bb \mdo
\vskip 0.0cm \noindent
\hhhh \mchoose the first such $f$
\vskip 0.0cm \noindent
\hhhh $h:=h-\left(a/lc(f)\right)\partial_\beta f$ where $\partial_\beta ld(f)=\theta$
\vskip 0.0cm \noindent
\hhh \mend
\vskip 0.0cm \noindent
\hh \mend
\vskip 0.3cm

\noindent
Correctness and termination of this algorithm is an obvious consequence of
Definition~\ref{inv_red} and correctness and termination of the polynomial
$L-$normal form algorithm~\cite{GB1}.

\begin{definition} A finite set $F$ is called {\em $L-$autoreduced} if
every $f\in F$ is irreducible modulo any other element $g\in F$.
An $L-$autoreduced set $F$ is called {\em $(L-)$involutive} if
$$
(\forall f\in F)\ (\forall \alpha\in \N^n) \ [\ NF_L(\partial_\alpha f,F)=0\ ].
$$
Given a derivative $\vartheta$ and an $L-$autoreduced set $F$, if there exist
$f\in F$ such that $ld(f)\prec_c \vartheta$ and
\begin{equation}
(\forall f\in F)\ (\forall \alpha \in \N^n)\ \left(\partial_\alpha ld(f)\prec_c
\vartheta\right)\ \
[\ NF_L(\partial_\alpha f,F)=0\ ]\,, \label{cond_pinv}
\end{equation}
then $F$ is called {\em partially involutive up to the derivative
$\vartheta$} with respect to the ranking $\prec_c$. $F$ is still said to
be partially involutive up to $\vartheta$ if $\vartheta\prec_c ld(f)$ for all
$f\in F$.
\label{inv_basis}
\end{definition}

\begin{corollary} If $F\subset \R$ is an $L-$involutive set, then
 every monomial set $U_j\in \M$ $(1\leq j\leq m)$ associated with
$ld_j(F)$ is $L-$involutive.
\label{cr_compl}
\end{corollary}

\begin{proof} It follows immediately from
 Definitions~\ref{id_noetherian} and \ref{inv_basis}.
\qed
\end{proof}

\begin{theorem}
 An $L-$autoreduced set $F\subset \R$ is involutive
 with respect to a continuous involutive division $L$ iff the
 following (local) involutivity conditions hold
$$
 (\forall f\in F)\ \left(\forall x_i\in NM_L(f,F)\right)\
 \ [\ NF_L(\partial_{x_i}\cdot f,F)=0\ ]\,.
$$
Correspondingly, partial involutivity~(\ref{cond_pinv}) holds
iff
$$
(\forall f\in F)\ (\forall x_i\in NM_L(f,F))\ (\partial_{x_i} \cdot
ld(f) \prec_c \vartheta)\ \ [\ NF_L(\partial_{x_i}\cdot f,F)=0\ ]\,.
$$
\label{th_inv_cond}
\end{theorem}

\begin{theorem} If $F\subset \R$ is an $L-$involutive basis of
 $[F]$, then it is also a differential \Gr basis.
\label{th_nf}
\end{theorem}

\noindent
The following theorem and corollary give an involutive analogue of
Buchberger chain
criterion \cite{Buch85} in application to linear differential bases.

\begin{theorem}
Let $F$ be a finite $L-$autoreduced set of linear differential polynomials
with respect to a continuous involutive division $L$,  and $NF_L(p,F)$ be
an algorithm of $L-$normal form. Then the following are equivalent:

\begin{enumerate}
\item $F$ is an $L-$involutive differential basis of $[F]$.
\item For all $g\in F, x\in NM_L(g,F)$ there is $f\in F$ satisfying
 $\partial_{x}\cdot ld(g)=\partial_L ld(f)$
 and a chain of elements in $F$ of the form
 $$f\equiv f_k,f_{k-1},\ldots,f_0,g_0,\ldots,g_{m-1},g_m\equiv g$$
 such that
 $$ NF_L\left(S_L(f_{i-1},f_i),F\right)=NF_L\left(S(f_0,g_0),F\right)=
 NF_L\left(S_L(g_{j-1},g_j),F\right)=0$$
 where $0\leq i\leq k$ and $0\leq j\leq m$, $S(f_0,g_0)$ is the conventional
 differential S-polynomial~\cite{Ollivier} and
 $S_L(f_i,f_j)=\partial_{x}\cdot f_i-\partial_L f_j$
 is its special form which occurs in involutive algorithms.
\end{enumerate}
\label{inv_chain}
\end{theorem}

\begin{corollary}
Let $F$ be a finite $L-$autoreduced set, and let $\partial_x\cdot g
$ be a nonmultiplicative prolongation of $g\in F$. If the
following holds
$$
(\forall h\in F)\ (\forall \,\partial_\alpha)\ \left(\,\partial_\alpha ld(h)\cdot u\prec_c
ld(g\cdot x)\,\right)\ \ [\ NF_L(h\cdot u,F)=0\ ]\,,
$$
$$
     (\exists f,f_0,g_0\in F)
\left[
\begin{array}{l}
 ld(f)=\partial_\beta ld(f_0)\,,\ ld(g)=\partial_\gamma ld(g_0) \\[0.1cm]
 \partial_x\cdot ld(g)=\partial_L ld(f)\,,\ lcm\left(ld(f_0),ld(g_0)\right)\prec_c
\partial_x\cdot ld(g)
 \\[0.1cm]
 NF_L\bigl(\partial_\beta\cdot f_0,F\bigl)=
 NF_L\bigl(\partial_\gamma\cdot g_0,F\bigl)=0
\end{array}
\right]\,,
$$
then the prolongation $\partial_x\cdot g$ may be discarded in the course of
an involutive algorithm.
\label{cor_criterion}
\end{corollary}

\section{Completion Algorithm}

The below given algorithm {\bf MinimalLinearInvolutiveBasis} is a differential
analogue of the polynomial algorithm {\bf MinimalInvolutiveBasis} of  paper~\cite{G98}.
In so doing, the conventional (non-involutive) autoreduction which is performed
in line 2 of the latter algorithm omitted, as this autoreduction is
optional~\cite{G98}.

Validity of the involutive chain criterion used in lines 11 and 23 is
provided by Theorem~\ref{inv_chain} and Corollary~\ref{cor_criterion}.
The proof of correctness and termination of the differential algorithm is
identical to the proof for its polynomial analogue~\cite{GB2,G98}. It
follows, that if the main ranking $\succ$ is orderly, then, given a generating
set of linear differential polynomials and a constructive involutive division,
algorithm {\bf MinimalLinearInvolutiveBasis} computes a minimal differential
basis whenever the latter exists. If the division is noetherian, the basis
is computed for any main ranking.

Though the output basis for a noetherian division does not depend on the
completion ranking, the proper choice of the latter may increase efficiency
of computation.

\begin{remark} If the algorithm {\bf MinimalLinearInvolutiveBasis}
takes a conventional
differential \Gr basis of the ideal $[F]$ as an input, then it produces the minimal
involutive differential basis just by enlargement of the input set with its
irreducible nonmultiplicative prolongations if any. This enlargement is done in
the lower {\bf while}-loop.
\label{rem_purecomp}
\end{remark}

\vskip 0.3cm
\setcounter{cc}{00}
\noindent
\h Algorithm {\bf MinimalLinearInvolutiveBasis}
\vskip 0.2cm
\noindent
\h {\bf Input:} $F$, $L$, $\succ$ (main ranking), $\succ_c$
 (completion ranking)
\vskip 0.0cm
\noindent
\h {\bf Output:} $G$, a minimal involutive basis of $[F]$
\vskip 0.0cm
\noindent
\h \mbegin
\hln
\hh \mchoose $g\in F$ with the lowest $ld(g)$ w.r.t. $\prec$
\hln
\hh $T:=\{(g,ld(g),\emptyset)\}$;\ \ $Q:=\emptyset$;\ \ $G:=\{g\}$
\hln
\hh \mfore $f\in F\setminus \{g\}$\bb \mdo
\hln
\hhh $Q:=Q\cup \{(f,ld(f),\emptyset )\}$
\hln
\hh \mrepeat
\hln
\hhh $h:=0$
\hln
\hhh \mwhile $Q\neq \emptyset$\bb \mand $h=0$\bb \mdo
\hln
\hhhh \mchoose $g$ in $(g,\theta,P)\in Q$ with the lowest $ld(g)$ w.r.t. $\prec$
\hln
\hhhh $Q:=Q\setminus \{(g,\theta,P)\}$
\hln
\hhhh \mif $Criterion(g,\theta,T)$ is false \bb \mthen $h:=NF_L(g,G)$
\hln
\hhh \mif $h\neq 0$\bb \mthen $G:=G\cup \{ h \}$
\hln
\hhhh \mif $ld(h)=ld(g)$\bb \mthen $T:=T\cup \{(h,\theta,P\cap NM_L(h,G))\}$
\hln
\hhhh \melse $T:=T\cup \{(h,ld(h),\emptyset )\}$
\hln
\hhhhh \mfore  $f$ in $(f,\vartheta,S)\in T$ s.t. $ld(f) \succ ld(h)$\bb \mdo
\hln
\hhhhhh $T:=T\setminus \{(f,\vartheta,S)\}$;\ \ $Q:=Q \cup \{(f,\vartheta,S)\}$;\ \
    $G:=G\setminus \{f\}$
\hln
\hhhhh \mfore  $(f,\vartheta,S)\in T$ \bb \mdo
\hln
\hhhhhh $T:=T\setminus \{(f,\vartheta,S)\}\cup  \{(f,\vartheta,S\cap NM_L(f,G))\}$
\hln
\hhh \mwhile exist $(g,\theta,P)\in T$ and $x\in NM_L(g,G)\setminus P$ and,
     if $Q \neq \emptyset$,
\hln
\hhhh s.t. $ld(\partial_x\cdot g) \prec ld(f)$ for all $f$ in
 $(f,\vartheta,S)\in Q$\bb \mdo
\hln
\hhhh \mchoose such $(g,\theta,P),x$ with the lowest $ld(\partial_x\cdot g)$
 w.r.t. $\prec_c$
\hln
\hhhhh $T:=T\setminus \{(g,\theta,P)\} \cup \{(g,\theta,P\cup \{x\})\}$
\hln
\hhhh \mif $Criterion(\partial_x\cdot g,\theta,T)$ is false \bb \mthen
 $h:=NF_L(\partial_x\cdot g,G)$
\hln
\hhhhh \mif $h\neq 0$\bb \mthen $G:=G\cup \{h\}$
\hln
\hhhhhh \mif $ld(h)=ld(\partial_x\cdot g)$ \bb \mthen
 $T:=T\cup \{(h,\theta,\emptyset )\}$
\hln
\hhhhhh \melse $T:=T\cup \{(h,ld(h),\emptyset )\}$
\hln
\hhhhhhh \mfore  $f$ in $(f,\vartheta,S)\in T$\bb with\bb
 $ld(f) \succ ld(h)$\bb \mdo
\hln
\hhhhhhhh $T:=T\setminus \{(f,\vartheta,S)\}$;\ \ $Q:=Q \cup
 \{(f,\vartheta,S\})$;\ \ $G:=G\setminus \{f\}$
\hln
\hhhhhhh \mfore  $(f,\vartheta,S)\in T$ \bb \mdo
\hln
\hhhhhhhh $T:=T\setminus \{(f,\vartheta,S)\}\cup  \{(f,\vartheta,S\cap NM_L(f,G))\}$
\hln
\hh \muntil $Q\neq \emptyset$
\hln
\h \mend
\hln
\vskip 0.2cm

\noindent
$Criterion(g,\theta,T)$ is true if there is $(f,\vartheta,S)\in T$
such that $ld(g)=\partial_L ld(f)$ and $lcm(\theta,\vartheta) \prec_c ld(g)$.

\begin{example}~\cite{Janet} The well-known Janet example with three independent and one
dependent variables $(n=3,m=1)$:
$$
\left\{
\begin{array}{l}
\partial_{11} y-x_2\partial_{33} y=0\,, \\
\partial_{22} y=0\,.
\end{array}
\right.
$$
The above completion algorithm applied for Janet, Pommaret and Lexicographically induced
divisions gives the following involutive bases, which coincide for both
pure lexicographical and graded lexicographical main rankings compatible
with~(\ref{var_order}) and which sorted in the descending
lexicographical order:

\begin{center}
\bg {|c|c|c|c|} \hline\hline
\Gr   & \multicolumn{2}{c|}{Involutive Bases}
\\ \cline{2-3}
basis &  Janet \& Pommaret & Lex. Induced \\ \hline
$\partial_{11} y-x_2\partial_{33} y$ & $\partial_{11} y-x_2\partial_{33} y$
&$\partial_{112}y-\partial_{33}y$ \\
$\partial_{22}y$ & $\partial_{122}y$ & $\partial_{11333}y$ \\
$\partial_{233}y$ & $\partial_{1233}y$ & $\partial_{1133}y$ \\
$\partial_{3333}y$ & $\partial_{13333}y$ &$\partial_{113}y-x_2\partial_{333}y$ \\
& $\partial_{22}y$ & $\partial_{11}y-x_2\partial_{33}y$ \\
& $\partial_{233}y$ & $\partial_{223}y$ \\
& $\partial_{3333}y$ & $\partial_{22}y$ \\
& & $\partial_{2333}y$ \\
& & $\partial_{233}y$ \\
& & $\partial_{3333}y$ \\
\hline \hline
\nd
\end{center}

\noindent
The first column contains the reduced differential \Gr basis, and Janet and
Pommaret bases are identical for this example.
\label{Ex_Janet}
\end{example}

\section{Initial Value Problem}

The results of this section generalize to arbitrary $L-$involutive
linear systems those obtained in Riquier-Janet theory~\cite{Riquier,Janet,Thomas}, for
Janet and Thomas divisions,
as well as in the formal theory~\cite{Pommaret,Pommaret94} for Pommaret
division, on posing
an initial value problem providing uniqueness and existence of solutions.

\begin{definition}~\cite{Riquier,Janet} If $\theta \in ld(F)$ is a
leading derivative in $F\subset \R$, then $\partial_\alpha \theta$ is called a
{\em principal derivative}. A derivative
which is not principal is called {\em parametric}. The monomial set
$W=\{\cup_{j=1}^m W_j\ |\ W_j\subset \M\}$ associated by~(\ref{d-m})
with the set of
parametric derivatives is called a {\em complementary set} of $F$.
\label{def_pripar}
\end{definition}

\begin{proposition}
Given a ranking $\prec$, if set $F$ is a linear $L-$involutive basis of differential
ideal $[F]$,
then the sets of principal and parametric derivatives
$($complementary set$)$
related to $F$ depend only on $[F]$ and $\prec$ and
do not depend on the choice of involutive division $L$.
\end{proposition}

\begin{proof} It follows immediately from the fact that any involutive basis
is a \Gr basis (Theorem~\ref{th_nf}).
\qed
\end{proof}

\begin{lemma} $($decomposition lemma$)$ Given a noetherian division $L$
 and $L-$ involutive set $F\subset \R$, every subset $W_j$ in the
 complementary monomial set of $F$ related to $j$-th differential
 indeterminate $y_j$ $(1\leq j\leq m)$ can be decomposed as a disjoint
 union
\begin{equation}
W_j=\cup_{v\in V_j} vL_v, \quad L_v\subseteq L(v,U_j\cup\{v\}),
\label{decomp}
\end{equation}
where $U_j$ is the $L-$involutive monomial set (not necessarily
nonempty) associated with $ld_j(F)$, and $V_j\in \M$ is a finite subset.
\label{lemma}
\end{lemma}

\begin{proof} Let $U$ and $W$ be a pair of monomial sets
associated with the principal and parametric derivatives of a
differential indeterminate in $F$. The complementary
set $W$ can be written as a disjoint union~\cite{CLO}
\begin{equation}
W=W_0\cup W_1\cup \cdots \cup W_d \label{union_1}
\end{equation}
where $d$ is the dimension of monomial ideal $(U)$, $W_0$
is a finite set, and every $W_i$ $(1\leq i\leq d)$ is a finite
disjoint union\footnote{The union in~(\ref{union_2}) considered
in~\cite{CLO} is not necessarily disjoint. However, unions
in (\ref{union_1}) and (\ref{union_2}) apparently can be rewritten as
disjoint by appropriate choice of $W_0$ and components of $W_r$.}
\begin{equation}
W_r=W_{r_1}\cup W_{r_2}\cup \cdots \cup W_{r_k} \label{union_2}
\end{equation}
with
\begin{equation}
W_{r_s}=\{\ w_{r_s}x^{\alpha_1}_{i_{s_1}}\cdots x^{\alpha_r}_{i_{s_r}}\
|\ \alpha_t\in \N\,,\
1\leq t\leq r\ \}\quad (1\leq s\leq k). \label{w_i}
\end{equation}

\noindent
For every $v\in W_0$ we shall take $L_v=\{1\}$ in~(\ref{decomp}). Thus,
for $d=0$  the decomposition~(\ref{decomp}) $W=W_0=\cup_{v\in W_0} \{v\}$
holds trivially. If $d>0$ we consider the finite set
\begin{equation}
V=W_0\cup_{r=1}^d \cup_{s=1}^{k} \{w_{r_s}\}, \label{gen}
\end{equation}
where monomials $w_{r_s}$ generate $W_{r_s}$ in accordance with~(\ref{w_i}).

We claim that elements in set~(\ref{gen}), and the
decompositions~(\ref{union_1}), (\ref{union_2}) they determine can be written
such that the union in
$W=\cup_{v\in V}vL_v$ with $L_v\subseteq L(v,U\cup \{v\})$
is disjoint in accordance with~(\ref{decomp}). To prove the claim
we define the degree $q$ of set $U$ as
$q=\max\{deg(u)\ |\ u\in U\}$, and choose all the monomials $w_{r_s}$
generating $W_{r_s}$ in~(\ref{w_i}) such that $deg(w_{r_s})=q$. Obviously this
can always be done by appropriate choice of $W_0$. Let now $V_1$ be
the set $V_1=\cup_{r=1}^d \cup_{s=1}^k \{w_{r_s}\}$, and let
$\hat{U}$ be a finite $L-$autoreduced completion of $U\cup V$. The existence
of $\hat{U}$ is guaranteed by noetherity of $L$. Now consider the set
$\hat{V}=\hat{U}\cap W \supseteq V_1$.
Its $L-$involutivity and property (d) of $L$ in Definition~\ref{ind_div} imply
$$ (\forall w\in W\setminus W_0)\ (\exists v\in \hat{V})\ \
[\ w\in vL(v,T)\subseteq vL(v,U\cup \{v\})\ ]. $$

\noindent
Thus, we obtain the desired decomposition $W=W_0\cup_{v\in \hat{V}} vL(v,\hat{U})$.
Disjointedness of this union follows from that
in~(\ref{union_2}) and Definition~\ref{autored_ms} of $L-$ autoreduction.
This proves the claim and the lemma.
\qed
\end{proof}

\begin{definition} Those elements $v_{j_k}$ (parametric derivatives)
which, in accordance with (\ref{decomp}),
generate the whole complementary set $W$, will be called {\em $L-$generators}
of the set.
The multiplicative variables $x_i$ satisfying $x_i\in L_{j_k}$
will be called {\em $(L-)$multipliers} of the generator $v_{j_k}$ and the
remaining variables will be called its {\em $(L-)$nonmultipliers}. The whole
set of $L-$generators of $W$ will be denoted by $\mathcal{G}_L$, and in accordance
with~(\ref{decomp})
\begin{equation}
{\mathcal{G}_L}=\cup_{j=1}^m V_j \label{full_gen}\,.
\end{equation}
\label{def_gen}
\end{definition}

\noindent
For a non-noetherian division $L$ a complementary set may not have a finite
set of $L-$generators as the following example shows.

\begin{example} Let involutive division $L$ be defined on $\M$ as follows.
 Variables $x_1,\ldots,x_{n-1}$ are separated into multiplicative and
 nonmultiplicative by Definition~\ref{div_J}. Let the variable $x_n$
 be also separated by Definition~\ref{div_J} if $deg_n(u)=0$ and $u\neq 1$,
 whereas if
 $deg_n(u)>0$ or if $u=1$, $x_n$ be nonmultiplicative for $u$.
 Then, the monomial set $U=\{x_1^2,x_1x_2,x_2\}$ is $L-$involutive in
 $K[x_1,x_2,x_3]$. Its complementary set has the infinite set of
 $L-$generators: $\mathcal{G}_L=\{1\}\cup \{x_1\}\cup_{i=1}^\infty \{x_3^i\}$.
\end{example}

\begin{remark} Decomposition~(\ref{decomp}) and the underlying
$L-$generator set~(\ref{def_gen}) are not uniquely defined, and usually
a more compact set $\mathcal{G}$ of $L-$generators (with less number of elements) than
that constructed in the proof of Lemma~\ref{lemma} can be chosen.
For example, for a Janet basis, $\mathcal{G}_P$ can always be
chosen~\cite{Janet} as union~(\ref{full_gen}) of sets $V_j$ such that
\begin{equation}
 (\forall \,V_j)\ (\forall v\in V_j)\ [\ L_v=J(v,U_j\cup \{v\}\ ]\quad
 (1\leq j\leq m), \label{J_gen}
\end{equation}
where $J$ stands for the Janet set of multiplicative monomials. Since for
$\hat{U}_j$,
as it constructed in the proof, the inclusion
$U_j\cup \{v\} \subseteq \hat{U}_j$
holds, the property (d) in Definition~\ref{inv_div} implies
$J(v,\hat{U}_j)\subset J(v,U_j\cup \{v\})$. Therefore, the set of
Janet generators
defined by~(\ref{J_gen}) is a subset of that constructed
in the proof of Lemma~\ref{lemma}.

For a Pommaret basis in the formal theory~\cite{Pommaret}
decomposition~(\ref{decomp})
is taken in the form
\begin{equation}
W=W_0\cup_{\{v\in W\ |\ deq(v)=q\}}vP(v), \label{P_decomp}
\end{equation}
where $P(v)$ denotes the set of Pommaret multiplicative monomials
for $v$, and $q$, as
in the proof, is the degree of the basis. The number
of Pommaret generators in~(\ref{P_decomp}) with $i$ multipliers
is called the $i$th
{\em Cartan character}\footnote{Cartan introduced these numbers in
his analysis of exterior PDEs~\cite{Cartan} and called them
{\em characters}.} $(1\leq i\leq n)$ of the basis and will be denoted by
 $\sigma_q^i$.
\label{rem_decomp}
\end{remark}

\begin{example} The complementary set of the monomial ideal $(U)$
 for $U=\{x_1^2x_3,x_1x_2,x_1x_3^2\}$ in Example~\ref{exm_1} is
$ W=\cup \{x_1^{i+1}\ |\ i\in \N\}\cup \{x_2^jx_3^k\ |\ j,k\in \N\}.$
Its most compact sets $\mathcal{G}_J$ and $\mathcal{G}_{D_{Lex}}$
together with their multipliers are:

\begin{center}
\bg {|c|c|c|c|} \hline\hline
\multicolumn{2}{|c|}{Janet division} &\multicolumn{2}{c|}
{Lex. induced division}
\\ \hline
Generator & Multipliers & Generator   & Multipliers \\ \hline
$1$      & $x_2,x_3$&  $1$   & $x_2,x_3$ \\

$x_1$    &   $-$   &  $x_1$ & $x_1$     \\
$x_1^2$  & $x_1$    &  $-$  & $-$     \\
\hline \hline
\nd
\end{center}

\noindent
We note that if the involutive bases $\tilde{U}_P,\tilde{U}_{D_{Lex}}$
given in Example~\ref{exm_1} are
sequentially enlarged with every single generator, then the sets of
Janet multipliers, in accordance with Remark~\ref{rem_decomp}, coincide
with the sets of multiplicative variables
$$ M_J(1,\tilde{U}_P\cup \{1\})=\{x_2,x_3\},\
M_J(x_1,\tilde{U}_P\cup \{1\})=\emptyset,\
M_J(x_1^2,\tilde{U}_P\cup \{x_1^2\})=\{x_1\}$$
whereas for lexicographically induced division, every
set of multipliers is the proper subset of multiplicative variables
$$ M_{D_{Lex}}(1,\tilde{U}_{D_{Lex}}\cup \{1\})=
M_{D_{Lex}}(x_1,\tilde{U}_{D_Lex}\cup \{x_1\})=\{x_1,x_2,x_3\}.$$
\end{example}

\begin{theorem} $($uniqueness theorem$)$ Let
${\cal F}$ be an $L-$involutive system of linear PDEs for an orderly
ranking. Then ${\cal F}$ has at most one solution satisfying the
following initial conditions: the derivatives associated with
$L-$generators of the complementary monomials are arbitrary functions
of their multipliers at the fixed values of their
nonmultipliers from coordinates of the initial point
$x_i=x^{\rm o}_i$ $(1\leq i\leq n)$), whereas the generators without
multipliers are considered to be arbitrary constants.
\label{th_uniqueness}
\end{theorem}

\begin{proof} Involutivity of ${\cal F}$ with respect to an orderly
ranking implies that the associated complementary monomial set
contains all the
monomials associated with the parametric derivatives. This statement is an
immediate consequence of the well-known fact~\cite{CLO} that for a
graded monomial ordering the Hilbert function of a polynomial ideal
is defined by the monomial ideal generated by the leading
monomials of a \Gr basis of the polynomial ideal.

Furthermore, by association~(\ref{d-m}), the decomposition~(\ref{union_1})
yields that every parametric derivative associated with a monomial
in $W\setminus W_0$ is produced by differentiation of
the uniquely defined parametric derivative ($L-$generator)
with respect to its multipliers. Assigning the fixed values
to all these parametric derivatives is obviously equivalent to
fixing some function of the multipliers. Therefore, given initial point $x_i=x^{\rm o}_i$,
in addition to
the set of arbitrary constants which associated with elements in $W_0$,
all the
parametric arbitrariness is determined by functions corresponding to the $L-$generators
and which are
arbitrary functions of the multipliers at the fixed values of
nonmultipliers from coordinates of the initial point.
\qed
\end{proof}

\begin{remark} Different involutive divisions give obviously
equivalent forms of initial value problem providing the uniqueness
of solutions. However, given a system of PDEs with an infinite
set of parametric derivatives, the writingf of such initial conditions
in accordance with Theorem~\ref{th_uniqueness} may be more compact for one
division than for another. We demonstrate this fact by
examples given below.
\label{rem_uniqueness}
\end{remark}

\begin{theorem}$($existence theorem$)$ Let ${\cal F}$ be an $L-$involutive linear
system for an orderly ranking, and let its coefficients be analytic
functions in an initial point $(x_i=x^{\rm o}_i)$. Then ${\cal F}$ has
precisely one solution which is analytic in this point if all arbitrary
functions in the initial data specified in Theorem~\ref{th_uniqueness} are
analytic in their arguments taking values from coordinates of the
initial point.
\label{th_existence}
\end{theorem}

\begin{proof} This is identical to the existence proof in
Riquier-Janet theory~\cite{Riquier,Janet} (see also~\cite{Ritt}).
\qed
\end{proof}

\begin{example} The complementary monomial set for Janet system in
Example~\ref{Ex_Janet} is finite and consists of 12 elements
$$W=\{1,x_1,x_2,x_3,x_1x_2,x_1x_3,x_2x_3,x_3^2,x_1x_2x_3,x_1x_3^2,x_3^3,x_1x_3^3\}.$$
By Theorem~\ref{th_existence}, its general solution depends on 12
arbitrary constants.
\label{Ex_Janet_1}
\end{example}

\begin{example}~\cite{Pommaret} The system of the first order PDEs
with four independent and one dependent variables $(n=4,m=1)$ and its
completion to involution for Janet or Pommaret division for any
ranking compatible with~(\ref{var_order}) are given by
$$
\left\{
\begin{array}{l}
\partial_1y+x_2\partial_3y+y=0,\\
\partial_2y+x_1\partial_4y=0,
\end{array}
\right.
\quad {\buildrel {J,P-}{\rm completion} \over {\buildrel \phantom{~~~~}
\over \Longrightarrow }} \quad
\left\{
\begin{array}{l}
\partial_1y+x_2\partial_4y+y=0,\\
\partial_2y+x_1\partial_4y=0,\\
\partial_3y-\partial_4y=0.
\end{array}
\right.
$$
The parametric derivatives $\partial_4^iy$ $(i\in \N)$ have
the only Janet generator $y\Longleftrightarrow 1$ with the only
multiplier $x_4$. Hence, the initial data providing the unique analytic
solution are
$y|_{x_1=x_1^{\rm o},x_2=x_2^{\rm o},x_3=x_3^{\rm o}}=\phi(x_4)$
with arbitrary function $\phi(x_4)$, analytic at $x_4=x_4^{\rm o}$.
The system can explicitly be integrated and its general solution is
$$y=e^{-(x_4-x_4^{\rm o})}\phi(x_3+x_4-x_1x_2+x_1^{\rm o}x_2^{\rm o}-
 x_3^{\rm o}). $$
The Pommaret generators are $y$ and $\partial_4 y$ without multipliers and
with multiplier $x_4$, respectively. This leads to the initial value problem
$$ y|_{x_1=x_1^{\rm o},x_2=x_2^{\rm o},x_3=x_3^{\rm o}}=c,\quad
\partial_4 y|_{x_1=x_1^{\rm o},x_2=x_2^{\rm o},x_3=x_3^{\rm o}}=\psi(x_4)$$
with arbitrary constant $c$ and arbitrary function $\psi$.
This shows that the Janet initial conditions are written in a more compact
form than those of Pommaret.
\label{Ex_Pommaret}
\end{example}

\begin{example}~\cite{Lewy} The well-known Lewy example with
 $n=3, m=2$ and $\eta_1,\eta_2\in \K$
$$
\left\{
\begin{array}{l}
\partial_1 y_1 - 2x_3\,\partial_2 y_1 - \partial_3 y_2 -
 2x_1\,\partial_2 y_2 = \eta_1(x_1,x_2,x_3),\\
\partial_1 y_2 + 2x_1\,\partial_1 y_1 + \partial_3 y_1 -
 2x_3\,\partial_2 y_2 = \eta_2(x_1,x_2,x_3).
\end{array}
\right.
$$
This system is involutive for any of Janet, Pommaret or lexicographically
induced divisions
and the orderly ranking with $\partial_1 y_j\succ \partial_2 y_j
\succ \partial_3 y_j$, $y_1\succ y_2$. Janet generators
are $y_1,y_2$. Each of them has multipliers $x_2,x_3$. This implies
the initial data
providing the uniqueness:
$y_j |_{x_1=x_1^o}=\phi_j(x_2,x_3)$ $(j=1,2)$
with arbitrary functions $\phi_j(x_2,x_3)$. Pommaret and lexicographically
induced divisions lead to a less compact writing of these conditions.
\label{Ex_Lewy}
\end{example}

\begin{remark} As shown by Lewy~\cite{Lewy} for Example~\ref{Ex_Lewy},
there exist the $C^\infty$
functions $\eta_1,\eta_2$ such that the system has no $C^\infty$ (and even $C^1$)
solutions. Therefore, analyticity in the Theorem~\ref{th_existence} statement
can not be replaced by smoothness.
\label{rem_analyticity}
\end{remark}

\noindent
We conclude this section with explicit formulae for the Hilbert function $HF_{[F]}$ and
Hilbert polynomial $HP_{[F]}$ of differential ideal $[F]$ represented by its linear
involutive basis $F$. These formulae are valid
for any involutive division and an orderly ranking.
For ordinary differential ideals, that is, for the case of single differential
indeterminate $(m=1)$, by association~(\ref{d-m}), they are the same as in
commutative algebra~\cite{Apel,GBC}. For partial differential case they involve
the number $m$ of differential indeterminates
\begin{eqnarray}
&&HF_{[F]}(s)=m\left(\begin{array}{c} n+s \\ s\end{array}\right) -
 \sum_{j=1}^m \sum_{i=0}^s \sum_{u\in U_j} \left(\begin{array}{c} i-deg(u)+\mu(u)-1 \\
\mu(u)-1 \end{array}\right), \label{HF}\\
&&HP_{[F]}(s)=m\left(\begin{array}{c} n+s \\ s\end{array}\right) -
 \sum_{j=1}^m \sum_{u\in U_j} \left(\begin{array}{c} s-deg(u)+\mu(u) \\
\mu(u) \end{array}\right). \label{HP}
\end{eqnarray}
Here $n$ is the number of independent variables, $U_j$ is the
monomial set associated with the set of leading derivatives $ld_j(F)$,
and $\mu(u)$ is the number of multiplicative elements of $u$.

The first term in the right hand side of~(\ref{HF}) is the
total number of derivatives of order $\leq s$. The triple
sum counts the number of principal derivatives among them
in accordance with Definition~\ref{inv_basis} which says that any
principal derivative is uniquely obtained by the multiplicative
prolongation of one of the leading derivatives in $F$. Thus,
(\ref{HF}) gives the number of parametric derivatives of
order $\leq s$, and for $s$ large enough it becomes polynomial~(\ref{HP}).

In the formal theory~\cite{Pommaret,Pommaret94} the Janet formula
is used:
$$ HP_{[F]}=\sum_{i=1}^n \left(\begin{array}{c} s-q+i-1 \\
i-1 \end{array}\right) \sigma_q^i. $$
Here the Hilbert polynomial~\cite{Janet}
is written in terms of Cartan characters $\sigma_q^i$
(see Remark~\ref{rem_decomp}). Apparently, this is~(\ref{HP}), rewritten
for Pommaret division in terms of Cartan characters.

\section{Lie Symmetry Analysis of PDEs}

Lie symmetry methods and their computerization yield a powerful practical tool for
analysis of nonlinear differential equations (see the review article~\cite{Hereman} and
references therein for more details). We present here the basic
computational formulae
and demonstrate, by two simple examples with a single nonlinear
evolution equation,
application of the above described involutive methods to finding the
classical infinitesimal symmetries.

\noindent
Given a finite system of polynomial-nonlinear PDEs
\begin{equation}
f_k(x_i,y_j,\ldots,\partial_\alpha y_j)=0,\  \ (1\leq i\leq n,\ 1\leq j\leq m,\ 1\leq k\leq r)
\label{system_PDEs}
\end{equation}
one looks for one-parameter infinitesimal transformations
\begin{equation}
\left\{
\begin{array}{l}
\tilde{x}_i(\lambda)=x_i+\xi_i(x_i,y_j)\lambda+O(\lambda^2), \\
\tilde{y}_j(\lambda)=y_j+\eta_j(x_i,y_j)\lambda+O(\lambda^2),
\end{array}
\right. \quad (1\leq i\leq n,\ 1\leq j\leq m). \label{trans}
\end{equation}
The conditions of invariance of~(\ref{system_PDEs}) under
transformations~(\ref{trans}) are
\begin{eqnarray}
&& \hat{Z}^{(\alpha)}f_k(x_i,y_j,\ldots,\partial_\alpha y_j)|_{f_s=0}=0,
\quad (1\leq k,s\leq r) \label{inv_cond} \\
&& \hat{Z}^{(\alpha)}=\xi_i\partial_{x_i} +
\eta_j\partial_{y_j}+
\zeta_{j;i}\partial_{y_{j;i}} +
\cdots +
\zeta_{j;\alpha}\partial_{y_{j;\alpha}}, \label{Z_operator}
\end{eqnarray}
where $\partial_iy_j$ denoted by $y_{j;i}$, etc.\footnote{In this section
the summation over repeated indices is always assumed.} Functions
$\zeta_{j;\ldots}$ involved in the differential operator~(\ref{Z_operator}) are
uniquely computed in terms of functions
$\xi_i,\eta_j$ and their derivatives by means of the recurrence relations
$$
\begin{array}{l}
\zeta_{j;i}=D_i(\eta_j)-y_{j;q}D_i(\xi_q), \\[0.2cm]
\zeta_{j;i_1\ldots i_p}=D_{i_p}(\zeta_{j;i_1\ldots i_{p-1}})
   -y_{j;i_1\ldots i_{p-1}q}D_{i_p}(\xi_q)\,,
\end{array}
$$
where $D_i$ is the total derivative operator with respect to $x_i$
$$
 D_i=\partial_{i} + y_{j;i}\partial_{y_j}
   + y_{j;ik}\partial_{y_{j;k}} + \cdots
$$
The invariance conditions~(\ref{inv_cond}) produce the overdetermined
system of linear homogeneous PDEs
in $\xi_i,\eta_j$ which is called the {\em determining system}. Its
particular solution yields an infinitesimal operator of the symmetry group
\begin{equation}
\hat{Z}=\xi_i\,\partial_{x_i} +
\eta_j\,\partial_{y_j}, \label{infin_oper}
\end{equation}
and the general solution yields all the infinitesimal operators.

Given initial system~(\ref{system_PDEs}), integration of the determining
system is generally a bottleneck of the whole procedure of
constructing these symmetry operators, and completion the
system to involution is the most universal algorithmic method of its
integration~\cite{Hereman}.

\begin{example}~\cite{Gerdt97} Diffusion type equation $y_t+yy_x-ty_{xx}=0$ $(n=2,m=1)$.
The symmetry operator~(\ref{Z_operator}) of the form
\begin{equation}
\hat{Z}=\xi_1\,\partial_t +
\xi_2\,\partial_x+\eta\,\partial_y \label{symm_opp}
\end{equation}
satisfies the determining system
\begin{eqnarray*}
&& \partial_{yy}\xi_1=0,\ \ \partial_{yy}\xi_2=0,\ \
   t\,\partial_{yy}\eta-2\,t\,\partial_{xy}\xi_2-
   2\,y\,\partial_y\xi_2=0,\\
&& \partial_y\xi_1=0,\ \
 2\,t^2\partial_{xy}\eta-t^2\partial_{xx}\xi_2-y\,t\,\partial_x\xi_2 +
 t\,\partial_t\xi_2+
 {y}\,\xi_1-t\,\eta=0,\\
&& t\,\partial_{xx}\eta-y\,\partial_x\eta-\partial_t\eta=0,\ \
t^2\,\partial_{xx}\xi_1-y\,t\,\partial_x\xi_1+2\,t\,\partial_x\xi_2-
t\,\partial_t\xi_1-
\xi_1=0,\\
&& t\,\partial_{xy}\xi_1+\,\partial_y\xi_2=0,\ \
\partial_x\xi_1=0.
\end{eqnarray*}
By choosing the orderly degree-reverse-lexicographical ranking with
$\partial_y \succ \partial_x\succ \partial_t$, $\xi_1\succ \xi_2\succ \eta$
and applying the completion algorithm of Sect.~5,
we obtain the (Pommaret, Janet, lexicographically induced) involutive
system
\begin{eqnarray*}
 &&\partial_y\xi_1=0,\quad \partial_y\xi_2=0,\quad \partial_y\eta=0,\quad
 \partial_x\xi_1=0,\quad \partial_x\xi_2-\frac{1}{t}\,\xi_1=0,\\
 &&\partial_x\eta=0,\quad \partial_t\xi_1-\frac{1}{t}\,\xi_1=0,\quad
 \partial_t\xi_2-\eta=0,
  \quad \partial_t\eta=0.
\end{eqnarray*}
The generators of parametric derivatives
$\xi_1,\xi_2,\eta$ have no multipliers.  Hence, the general solution
depends on three arbitrary constants $c_1,c_2,c_3$, and it can easily be
obtained by explicit integration of the
involutive system
$$\xi_1=c_1t,\quad \xi_2=c_1x+ c_2t+c_3,\quad \eta=c_2.$$
Respectively, the Lie symmetry group is three-dimensional. Its
symmetry operators
$ \hat{Z}_1=t\partial_t+x\partial_x$,\ \
$\hat{Z}_2=t\partial_x+\partial_y$,\ \
$\hat{Z}_3=\partial_x$
form the Lie algebra
$[{\hat Z}_1,{\hat Z}_2]=0$,\ \ $[{\hat Z}_2,{\hat Z}_3]=0$,\ \
   $[{\hat Z}_1,{\hat Z}_3]=-{\hat Z}_3$.
\label{diff_eqn}
\end{example}

\begin{example}~\cite{AC} The Harry Dym equation $\partial_ty-y^3\partial_{xxx}y=0$ $(n=2,m=1)$
which was already used in~\cite{Hereman} as an illustrative example.
The symmetry operator in the form~(\ref{symm_opp}) is now determined by
the system
\begin{eqnarray*}
&& \partial_y\xi_1=0,\quad \partial_x\xi_1=0,\quad \partial_y\xi_2=0,
\quad \partial_{yy}\eta=0,\\
&& \partial_{xy}\eta-\partial_{xx}\xi_2=0,\quad
  \partial_t\eta - y^3\partial_{xxx}\eta=0,\\
&& 3\,y^3\partial_{xxy}\eta +\partial_t\xi_2-y^3\partial_{xxx}\xi_2=0,
 \quad y\,\partial_t\xi_1-3\,y\,\partial_x\xi_2+3\,\eta=0.
\end{eqnarray*}
Its Janet and Pommaret involutive form for the same ranking as in
the previous example is
\begin{eqnarray*}
&& \partial_{xx}\eta=0,\ \ \partial_{xt}\eta=0,\ \ \partial_y\eta-\frac{1}{y}\,{\eta}=0,\ \
\partial_t\eta=0,
\quad \partial_y\xi_2=0,\\
&& \partial_x\xi_2-\frac{1}{3}\,\partial_t\xi_1-
\frac{1}{y}\,\eta=0, \ \partial_t\xi_2=0,\ \ \partial_{tt}\xi_1=0,\quad \partial_y\xi_1=0,
\ \ \partial_x\xi_1=0.
\end{eqnarray*}
There are five generators of parametric derivatives
$\xi_1,\partial_t\xi_1,\xi_2,\eta,\partial_x\eta$ which
have no multipliers that implies the five-dimensional
Lie symmetry group. The involutive determining system in this example
is also easy to integrate:
$$\xi_1=c_1+c_2t,\quad \xi_2=c_3+c_4x+c_5x^2,\quad
\eta=(c_4-\frac{1}{3}\,c_2+2\,c_5x)\,y.$$
This gives the Lie symmetry operators
$$ Z_1=\partial_t,\ \ Z_2=t\,\partial_t-\frac{1}{3}y\,\partial_y,\
\ Z_3=\partial_x,\ \ Z_4=x\,\partial_x+y\,\partial_y,\ \
Z_5=x^2\partial_x+2xy\,\partial_y $$
with the following nonzero commutators of the symmetry algebra
$$ [Z_1,Z_2]=Z_1,\quad [Z_3,Z_4]=Z_4,\quad [Z_3,Z_5]=2\,Z_4,\quad
[Z_4,Z_5]=Z_5.$$  \label{HD_eqn}
\end{example}

\section{Conclusion}

Most of the above presented definitions, statements and constructive methods
can be extended to finite sets of differential polynomials in $\R$ which,
given a ranking, are linear with respect to their highest rank
(\,{\em principal}\,) derivatives. In Riquier-Janet theory the corresponding systems
of PDEs are
called {\em orthonomic}. Their completion to involution, for any constructive
and noetherian division, could be done much like linear systems. The
essential obstruction here is a non-orthonomic integrability condition.
Moreover, even if such an integrability condition is explicitly solvable
with respect to its principal derivative, then this leads to non-polynomial
orthonomicity, and, thereby, to difficulty in the use of constructive methods
of differential and commutative algebra. In the latter case some geometric
features of the formal theory may be useful for computational
purposes~\cite{Reid3}.

However, given an orthonomic system of polynomial PDEs and an involutive
division $L$, one can always verify if it is $L-$involutive. Analytic
involutive orthonomic systems admit posing an initial value problem
providing the existence and uniqueness of solution. One can, hence,
determine arbitrariness in the general solution as it is done in Sect.~6
for linear systems. In particular, the compact general formulae~(\ref{HF})
and~(\ref{HP}) for the Hilbert function and Hilbert polynomial are also
valid for involutive orthonomic equations.

We are going to implement the completion algorithm {\bf
MinimalLinearInvolutiveBasis} (\,Sect.~5\,) after examination and
optimization of its polynomial analogue~\cite{GB2,G98}. Though its
implementation in Reduce for Pommaret division~\cite{GB1} has already
shown its efficiency, the differential case needs more careful analysis
of implementation and optimization issues to be applicable to
PDEs of practical interest. Thus, in Lie symmetry analysis of relatively
small systems it is easy to obtain determining systems of many
hundreds equations. Currently, the most efficient completion
algorithm for linear systems implemented in some
packages for Lie symmetry analysis~\cite{Hereman} is that of
paper~\cite{Reid1}. Its underlying implementations allow to
treate hundreds determining equations (cf.~\cite{Reid3}).
As for significantly larger determining systems, they are hardly
tractable by the present day
computer algebra tools, whereas there are practical needs in it.
In gas dynamics, for instance, the group classification of
the system of five second order PDEs describing a viscous heat conducting
gas and involving five dependent and four independent variables (three spatial
and one temporal) \cite{Bublik}, leads to
the determining system containing more than 200 000 equations.

In our intention to extract, in the process of implementation,
the maximal possible efficiency from the algorithms proposed, we hope,
first of all, to detect (heuristically) the most optimal choice of
involutive division. As the first step in this direction an
implementation of the monomial completion for different divisions
has been done in {\em Mathematica} and for Janet division in C~\cite{GBCK}.


\begin{thebibliography}{99}

\bibitem{Cartan1} Cartan, E.: Sur certaines expressions diff\'{e}rentielles
\`{a} le probl\`{e}me de Pfaff. {\em Annales Ecole Normale}, 3-e serie,
{\bf 16}, 1899, 239-332; Sur l'integration des syst\`{e}mes d'\'{e}quations aux
diff\'{e}rentielles totales. Ibid., {\bf 18} (1901) 241-311.

\bibitem{Kahler} K\"{a}hler, E.: {\em Einf\"{u}hrung in die Theorie der Systeme
 von Differentialgleichungen}, Teubner, Leipzig, 1934.

\bibitem{Cartan} Cartan, E.: {\em Les Syst\`{e}mes Diff\'{e}rentielles
 Ext\'{e}rieurs et leurs Applications G\'{e}ometriques}, Paris, Hermann, 1945.

\bibitem{ASY} Arais, E.A., Shapeev, V.P., Yanenko, N.N.: Realization of
 Cartan's Method of Exterior Differential Forms on an Electronic Computer.
 {\em Sov. Math. Dokl.} {\bf 15(1)} (1974) 203-205.

\bibitem{HT} Hartley, D., Tucker, R.W.: Constructive Implementation of
 the Cartan-K\"{a}hler Theory of Exterior Differential Systems.
 {\em J. Symb. Comp.} {\bf 12} (1991) 655-667.

\bibitem{Riquier} Riquier, C.: {\em Les Syst\`emes d'Equations aux
D\'eriv\'ees Partielles}, Gauthier-Villars, Paris, 1910.

\bibitem{RR97} Rust, C.J., Reid G.J.: Rankings of Partial Derivatives.
 In: {\em Proceedings of ISSAC'97}, W.K\"{u}chlin (ed.), ACM Press, 1997
 pp. 9-16.

\bibitem{Janet} Janet, M.: {\em Le\c cons sur les Syst\`emes
 d'Equations aux D\'eriv\'ees Partielles}, Cahiers Scientifiques, IV,
 Gauthier-Villars, Paris, 1929.

\bibitem{Schwarz1} Schwarz, F.: The Riquier-Janet Theory and its
application to Nonlinear Evolution Equations. {\em Physica} {\bf
11D} (1984) 243-251.

\bibitem{Topunov} Topunov, V.L.: Reducing Systems of Linear
 Differential Equations to a Passive Form. {\em Acta Appl. Math.} {\bf 16}
 (1989) 191-206.

\bibitem{Schwarz2} Schwarz. F.: An Algorithm for Determining the Size
 of Symmetry Groups. {\em Computing} {\bf 49} (1992) 95-115.

\bibitem{Pommaret} Pommaret, J.F.: {\em Systems of Partial
 Differential Equations and Lie Pseudo\-groups}, Gordon \& Breach,
 New York, 1978.

\bibitem{Pommaret94} Pommaret, J.F.: {\em Partial Differential
 Equations and Group Theory. New Perspectives for Applications}, Kluwer,
Dordrecht, 1994.

\bibitem{Seiler1} Sch\"{u}, J., Seiler, W.M., Calmet, J.: Algorithmic
 Methods for Lie Pseudogroups. In: {\em Modern Group Analysis: Advanced
 Analytical and Computational Methods in Mathematical Physics},
 N.Ibragimov et al (eds.), Kluwer, Dordrecht, 1993, pp. 337-344.

\bibitem{Seiler2} Seiler, W.M.: {\em Applying AXIOM to Partial Differential
 Equations}. Internal Report 95-17, Universit\"{a}t Karlsruhe, Fakult\"{a}t
 Informatik, 1995.

\bibitem{Thomas} Thomas, J.: {\em Differential Systems}. AMS
 Publication, New York, 1937.

\bibitem{Boulier} Boulier, F., Lazard, D., Ollivier, F., Petitot M.:
Representation for the Radical of a Finitely Generated Differential
Ideal. In: {\em Proceedings of ISSAC'95},  A.H.M. Levelt (ed.),
ACM Press, 1995, pp. 158-166.

\bibitem{ZB93} Zharkov, A.Yu., Blinkov, Yu.A.: Involutive Approach
to Investigating Polynomial Systems. In: Proceedings of ``SC 93",
International IMACS Symposium on Symbolic Computation: New Trends
and Developments (Lille, June 14-17, 1993). {\em Math. Comp. Simul.} {\bf 42} (1996) 323-332.

\bibitem{GB1} Gerdt, V.P., Blinkov, Yu.A.: Involutive Bases of Polynomial
 Ideals. Preprint-Nr. 1/1996, Naturwissenschaftlich-Theoretisches
 Zentrum, University of Leipzig; {\em Math. Comp. Simul.} {\bf 45} (1998) 519-542.

\bibitem{Apel} Apel, J.: Theory of Involutive Divisions and an
Application to Hilbert Function. {\em J. Symb. Comp.} {\bf 25} (1998) 683-704.

\bibitem{GB2} Gerdt, V.P., Blinkov, Yu.A.: Minimal Involutive Bases.
 {\em Math. Comp. Simul.} {\bf 45} (1998) 543-560.

\bibitem{GBC} Gerdt, V.P., Berth, M., Czichowski, G.: Involutive
 Divisions in Mathematica: Implementation and Some Applications.
 In: {\em Proceedings of the 6th Rhein Workshop on Computer Algebra}
 (Sankt-Augustin, Germany, March 31 - April 3, 1998), J.Calmet (Ed.),
 Institute for Algorithms and Scientific Computing, GMD-SCAI,
 Sankt-Augustin, 1998, pp.74-91.

\bibitem{G98} Gerdt, V.P.: Involutive Division Technique:
 Some Generalizations and Optimizations, Preprint JINR E5-98-151,
 Dubna, 1998. To be published in the {\em Proceedings of "CASC'98"} (April
 20-24, 1998, St.Petersburg).

\bibitem{Reid1} Reid, G.J.: Algorithms for Reducing a System of PDEs to
 Standard Form, Determining the Dimension of its Solution Space and
 Calculating its Taylor Series Solution. {\em Euro. J. Appl. Maths.} {\bf
 2} (1991) 293-318.

\bibitem{Reid3} Reid, G.J., Wittkopf, A.D., Boulton A.: Reduction of
 Systems of Nonlinear Partial Differential Equations to Simplified
 Involutive Form. {\em Euro. J. Appl. Maths.} {\bf 7} (1996) 635-666.

\bibitem{Carra} Carra'Ferro, G.: Gr\"obner Bases and Differential
 Algebra. {\em Lec. Not. in Comp. Sci.} {\bf 356} (1987) 129-140.

\bibitem{Ollivier} Ollivier, F.: Standard Bases of Differential
 Ideals. {\em Lec. Not. in Comp. Sci.} {\bf 508} (1990) 304-321.

\bibitem{Hereman} Hereman, W.: Symbolic software for the computation
 of Lie
 symmetry analysis. In: {\em CRC Handbook of Lie Group Analysis of Differential
 Equations, Volume 3: New Trends in Theoretical Developments and Computational
 methods}, Ibragimov, N.H. et al. (eds.), CRC Press, Boca Raton, 1995, pp. 367-413.

\bibitem{Reid2} Reid, G.J.: Finding Abstract Lie Symmetry Algebras of
 Differential Equations without Integrating Determaning Equations.
 {\em Euro. J. Appl. Maths.} {\bf 2} (1991) 319-340.

\bibitem{Mansfield} Mansfield, E., Clarkson, P.A.: Application of the
 Differential Algebra Package {\tt diffgrob2} to Classical Symmetries
 of Differential Equations, {\em J. Symb. Comp.} {\bf 23} (1997) 517-533.

\bibitem{Ritt} Ritt, J.F.: {\em Differential Algebra}, AMS Publication, New York, 1950.

\bibitem{Kolchin} Kolchin, E.R.: {\em Differential Algebra and Algebraic
 Groups}, Academic Press, New York, 1973.

\bibitem{Gerdt97} Gerdt ,V.P.: Gr\"obner Bases and Involutive Methods for
 Algebraic and Differential Equations, {\em Math. Comp. Model.} {\bf 25}, No.8/9 (1997) 75-90.

\bibitem{GBCK} Gerdt, V.P., Berth, M., Czichowski, G., Kornyak V.V.:
 Construction of Involutive Monomial Sets for Different
 Involutive Divisions. {\em This volume}.

\bibitem{Buch85} Buchberger, B.: Gr\"obner Bases: an Algorithmic
 Method in Polynomial Ideal Theory. In: {\em Recent Trends in
 Multidimensional System Theory}, Bose, N.K. (ed.), Reidel, Dordrecht, 1985.
 pp. 184-232.

\bibitem{CLO} Cox, D., Little, J., O'Shea D.: {\em Ideals, Varieties and
 Algorithms}, 2nd Edition, Springer-Verlag, New York, 1996.

\bibitem{Lewy} Lewy. H.: An Example of a Smooth Linear Partial Differential
 Equation without Solution, {\em Ann. Math.} {\bf 66} (1957) 155-158.

\bibitem{AC}  Ablowitz, M.J., Clarkson, P.A.:
{\em Solitons, Nonlinear Evolution Equations and Inverse Scattering},
London Mathematical Society Lecture Notes on Mathematics {\bf 149},
Cambridge University Press, Cambridge, UK, 1991.

\bibitem{Bublik} Bublik, V.V.: Group Classification of Equations for Dynamics
of Viscous Heat Conducting Gas, In: {\em Dynamics of Continuous Medium}
{\bf 113}, Novosibirsk, 1998, pp. 19-21 (in Russian).

\end{thebibliography}
\end{document}